\documentclass[12pt]{amsart}
\usepackage{amssymb}
\usepackage{amsmath}
\usepackage{longtable}
\newcommand\g{{\mathfrak g}}
\newcommand\h{{\mathfrak h}}
\newcommand{\f}{\mathfrak{f}}
\newcommand\gl{\mathfrak{gl}}

\renewcommand\l{\mathfrak l}
\newcommand\z{\mathfrak z}
\newcommand\s{\mathfrak s}
\renewcommand{\t}{\mathfrak{t}}

\newcommand\tr{\operatorname{tr}}
\newcommand\C{\mathbb C}

\newcommand\N{\mathbb N}

\renewcommand\sl{\mathfrak{sl}}
\newcommand\so{\mathfrak{so}}
\renewcommand\sp{\mathfrak{sp}}
\newcommand{\ad}{\mathop{\rm ad}\nolimits}
\newcommand{\Ad}{\mathop{\rm Ad}\nolimits}
\newcommand{\Int}{\mathop{\rm Int}\nolimits}
\newcommand{\Aut}{\mathop{\rm Aut}\nolimits}

\newcommand{\rank}{\mathop{\rm rank}\nolimits}
\newcommand{\GL}{\mathop{\rm GL}\nolimits}
\newcommand\SL{\mathop{\rm SL}\nolimits}
\newcommand\Sp{\mathop{\rm Sp}\nolimits}
\newcommand\SO{\mathop{\rm SO}\nolimits}
\renewcommand\O{\mathop{\rm O}\nolimits}

\newcommand\Spec{\mathop{\rm Spec}\nolimits}
\newcommand\quo{/\!/}
\unitlength=1mm  \numberwithin{equation}{section}
\author{Ivan V. Losev}
\thanks{{\it 2000 Mathematics Subject Classification.} Primary 17B20, 14R20;
Secondary 14L30.}\thanks{{\it Key words and phrases.} Semisimple
Lie algebras, conjugacy of embeddings, invariants of sets of
elements in Lie algebras}
\title{On invariants of a set of elements of a semisimple Lie algebra}
\date{December 21,2005}
\newtheorem{Theorem}{Theorem}[section]
\newtheorem{Proposition}[Theorem]{Proposition}
\newtheorem{Corollary}[Theorem]{Corollary}
\newtheorem{Lemma}[Theorem]{Lemma}
\theoremstyle{definition}

\newtheorem{Definition}[Theorem]{Definition}
\newtheorem{Remark}[Theorem]{Remark}
\begin{document}
\begin{abstract}
Let $G$ be a complex reductive algebraic group, $\g$ its Lie
algebra and $\h$ a reductive subalgebra of $\g$, $n$ a positive
integer. Consider the diagonal actions $G:\g^n, N_G(\h):\h^n$. We
study a relation between the algebra $\C[\h^n]^{N_G(\h)}$ and its
subalgebra consisting of restrictions to $\h^n$ of elements of
$\C[\g^n]^G$.
\end{abstract}
\maketitle \tableofcontents
\section{Introduction}
Let $G$ be a  reductive algebraic group over the field $\C$, $\g$
its Lie algebra, $\h$ a reductive algebraic subalgebra of $\g$ and
$\widetilde{H}=N_G(\h)$ the normalizer of $\h$ in the group $G$.

Let $n$ be a positive integer. One has the diagonal actions
$G:\g^n, \widetilde{H}:\h^n$. Consider a subalgebra of $\C[\h^n]$
whose elements are restrictions of elements of $\C[\g^n]^G$ to
$\h^n$. We denote this algebra by $\C[\h^n]^G$. Clearly,
$\C[\h^n]^G\subset \C[\h^n]^{\widetilde{H}}$. It is interesting to
ask how different these two algebras can be.

It is more convenient to translate this question into geometric
language. As usual, we denote  by $\h^n\quo \widetilde{H}$ the
categorical quotient for the action $\widetilde{H}:\h^n$. In other
words, $\h^n\quo \widetilde{H}=\Spec(\C[\h^n]^{\widetilde{H}})$.
Put $\h^n\quo G=\Spec(\C[\h^n]^G)$. The inclusion
$\C[\h^n]^G\hookrightarrow \C[\h^n]^{\widetilde{H}}$ induces the
morphism of algebraic varieties $\psi_n:\h^n\quo H\rightarrow
\h^n\quo G$. The aim of this paper is to answer the following
questions: is $\psi_n$ isomorphism, birational, bijective, finite
morphism?

The property of  $\psi_n$ being bijective has a nice alternative
description. Namely, $\psi_n, n\geqslant 2,$ is bijective iff for
any reductive algebraic group $F$ and  any homomorphisms
$P_1,P_2:F\rightarrow H$ the following conditions are equivalent:
\begin{enumerate}
\item There is $g\in G$ such that $\Ad(g)\circ \rho_1=\rho_2$.
\item There is $g\in \widetilde{H}$ such that
$\Ad(g)\circ\rho_1=\rho_2$.
\end{enumerate}
Here $\rho_i=dP_i:\f\rightarrow \h,i=1,2,$ is the corresponding
homomorphism of Lie algebras.

The starting point for our work is E.B. Vinberg's
paper~\cite{Vinberg}, where  the morphism $\Psi_n:H^n\quo
\widetilde{H}\rightarrow H^n\quo G$ defined analogously to
$\psi_n$ was studied (here $H$ is a reductive subgroup of $G$).
The main result of that paper is that $\Psi_n$ is the morphism of
normalization for $n\geqslant 2$.

Now we list our main results.

At first, the morphism $\psi_n$ is always finite
(Proposition~\ref{prop.1.2}). If $n>1$ it is also birational
(Proposition~\ref{prop1.1}). Therefore, $\psi_n$ is the morphism
of normalization for $n>1$. In general, $\psi_1$ is not
birational. However, the following statement holds
\begin{Theorem}\label{Thm:0.1} Suppose that $G=\GL_n$ and $\h$ is a simple algebra different from
$\so_9,\sp_8, \so_{16}$, $\sl_8,\sl_9$. Then $\psi_1$ is
birational. If $\h$ is one of five algebras listed above, then for
some positive integer $m$ there exists an embedding
$\h\hookrightarrow \gl_m$ such that the corresponding morphism
$\psi_1$ is not birational.
\end{Theorem}

If $\psi_n$ is bijective, then $\psi_1$ is also bijective. For
some $G$ and $\h$ the converse is true. To describe such pairs we
need the following definition:

\begin{Definition}\label{defi:1}
Let $H$ be a reductive algebraic group and $\h$ be its Lie
algebra. Suppose that $\h$ does not contain simple ideals
isomorphic to $E_6,E_7,E_8$. Further, suppose that for every ideal
 $\h_1\subset \h$ isomorphic to $\so_{2k},k>3,$ there exists $h\in H$ such
that the restriction of $\Ad(h)$ to a simple ideal $\h_2\subset
\h$ is an outer involutory automorphism of $\h_2$ if $\h_2=\h_1$
and the identity otherwise (the claim of $\Ad(h)|_{\h_1}$ being
involution is essential only for $\h_1\cong \so_8$). Then we say
that $H$ is a {\it group of type I}.
\end{Definition}

Put $\overline{H}=\widetilde{H}/Z_G(\h)$.
\begin{Theorem}\label{Thm1:2}
If $\overline{H}$ is a group of type I and $\psi_1$ is bijective,
then so is $\psi_n$.
\end{Theorem}

However, in some cases $\psi_2$ is not bijective.

\begin{Proposition}\label{Thm1:3}
If $G$ is a group of type I, $\h$ is a simple algebra and
$\overline{H}$ is not a group of type I, then $\psi_n$ is not
bijective for $n>1$.
\end{Proposition}

In fact, the claim of $\h$ being simple can be omitted, but we do
not prove it.

Suppose now that $G=\GL_m$. It is known from classical invariant
theory that the algebra $\C[\g^n]^G$ is generated by polynomials
of the form $f(X_1,\ldots, X_n)=\tr(X_{i_1}X_{i_2}\ldots
X_{i_k})$. Clearly, if $\psi_n$ is isomorphism, then $\psi_1$ is
isomorphism too. In some cases the opposite is true.
\begin{Proposition}\label{prop:0.1} The algebra $\C[\h^n]^G$ is generated by elements
of the form $f(L(X_1,\ldots,X_n))$, where $f\in\C[\h]^G$ and $L$
is a Lie polynomial on $X_1,\ldots, X_n$.
\end{Proposition}

Using Proposition~\ref{prop:0.1} one can prove that $\psi_n$ is
isomorphism if so is $\psi_1$ for $\h=\sp_{2m}$, $\h=\so_{2m+1}$,
$\h=\so_{2m}$ with $\overline{H}\cong\Ad(\O_{2m})$, $\h=\sl_m$
with $\overline{H}=\Ad(\SL_m)$.

The author is grateful to  E.B. Vinberg for constant attention to
this work and to A.N. Minchenko for useful discussions.

\section{Finiteness and birationality of $\psi_n$ in general case}\label{SECTION_1}
Let $G,\g,\h,\widetilde{H}$ be as above and $n$ be a positive
integer. The natural map $\Phi_n:\h^n\rightarrow \g^n\quo G$ is
constant on $\widetilde{H}$-orbits. Therefore there is a natural
morphism $\varphi_n:\h^n\quo \widetilde{H}\rightarrow \g^n\quo G$.
The closure of $\varphi_n(\h^n\quo H)$ coincides with $\h^n\quo
G$.

The following proposition is due to  Richardson~\cite{Richardson}.
\begin{Proposition}\label{prop1}
Let $G$ be a reductive group, $\g$  its Lie algebra and $n$ a
positive  integer. Consider the action $G:\g^n$ as above. The
orbit of an $n$-tuple ${\bf x}=(x_1,\ldots,x_n)\in\g^n$ is closed
(resp. contains 0 in its closure) iff the algebraic subalgebra of
$\g$ generated by $x_1,\ldots,x_n$ is reductive (resp. consists of
nilpotent elements).
\end{Proposition}
\begin{Corollary}\label{cor1}  $\h^n\quo \widetilde{H}$ (respectively, $\h^n\quo G$) can
be identified with a set of equivalence classes of $n$-tuples
$(x_1,\ldots, x_n)\in\h^n$ generating a reductive subalgebra of
$\h$ modulo $\widetilde{H}$- (respectively, $G$-) conjugacy.
\end{Corollary}

\begin{Lemma}\label{lem2}
Any reductive algebraic Lie algebra $\h$ can be generated by two
elements (as an {\it algebraic} algebra). If  $\h$ is commutative,
then it can be generated by one element.
\end{Lemma}
\begin{proof}
The proof is completely analogous to the proof of Proposition 2
in~\cite{Vinberg}.
\end{proof}

\begin{Proposition}\label{prop1.1}
Suppose $n>1$. Then  $\psi_n$ is birational.
\end{Proposition}
\begin{proof}
The proof  is completely analogous to the proof of Theorem 2
in~\cite{Vinberg}.
\end{proof}

\begin{Proposition}\label{prop.1.2}
The morphism $\psi_n$ is finite.
\end{Proposition}
\begin{proof}
Denote by $\mathcal{N}_G$ (resp. $\mathcal{N}_{\widetilde{H}}$)
the null-cone for the action $G:\g^n$ (resp.
$\widetilde{H}:\h^n$). It follows from Proposition~\ref{prop1},
that $\mathcal{N}_G\cap \h^n=\mathcal{N}_{\widetilde{H}}$. Now our
statement follows from one version of Nullstellensatz
(see~\cite{Kraft}, ch.2, s. 4.3, Theorem 8).
\end{proof}
\begin{Corollary}
Suppose $n>1$. Then $\psi_n$ is the morphism of normalization.
\end{Corollary}
\begin{proof}
Since $\h^n\quo \widetilde{H}$ is normal, this is a direct
consequence of Propositions~\ref{prop1.1},\ref{prop.1.2}.
\end{proof}

\section{Birationality of $\psi_1$ for $G=\GL_n$ and simple algebra
$\h$}\label{SECTION_2}

In this section we prove Theorem~\ref{Thm:0.1}. Here we assume
that $\h$ is a simple Lie algebra of rank $r$. Denote by
$\alpha_1,\ldots,\alpha_r$  simple roots of $\h$, by
$\pi_1,\ldots,\pi_r$ the corresponding fundamental weights and by
$P$ and $Q$ the weight and the root lattices of $\h$,
respectively.

\begin{Lemma}\label{lem3}
Let $\Delta$ be an irreducible root system in a real vector space
$V$, $W$  its Weyl group, $Q$ the lattice generated by $\Delta$.
Denote by $(\cdot,\cdot)$ $W$-invariant scalar product on $V$. Let
$g\in \Aut(Q)$ be an orthogonal linear operator. If $\Delta\neq
C_4$, then $g\in\Aut(\Delta)$.
\end{Lemma}
\begin{proof}
It follows from the construction of root systems (see, for
example,~\cite{Bourbaki}, chapter 6) that if $\Delta\neq C_l$,
then elements of $Q$ lying in $\Delta$ are all elements satisfying
some conditions on their length. For $\Delta=C_l$ the same is true
for roots of minimal length. Suppose now that $\Delta=C_l$ and
$l\neq 4$. Let $g$ be an element of $\mathrm{O}(V)$ such that the
set of elements of $\Delta$ of minimal length is invariant under
$g$. Then $g\in W$ and we are done.
\end{proof}

\begin{proof}[Proof of Theorem~\ref{Thm:0.1}]
Denote by $\t$ a Cartan subalgebra of $\h$. The points of
$\h\quo\widetilde{H}\cong\h\quo \overline{H}$ (respectively,
$\h\quo G$) are in one-to-one correspondence with equivalence
classes of semisimple elements of $\h$ modulo $\overline{H}$-
(respectively $G$-) conjugacy (see Corollary~\ref{cor1}). If $t$
is an element of $\t$ in general position, then $gt\in\t$ for some
$g\in G$ implies $g\in N_G(\t)$. Thus, $\psi_1$ is birational iff
$N_G(\t)/Z_G(\t)=N_{\overline{H}}(\t)/Z_{\overline{H}}(\t)$.

Suppose that $\h\neq \sp_8,\so_9,\sl_8,\sl_9,\so_{16}$. Denote by
$\varphi$ an embedding of $\h$ into $\gl_m$. We identify $\h$ and
$\varphi(\h)$. Denote by $N$ the image under the natural
homomorphism of $N_G(\t)$ in $\GL(\t)$. It is clear that $N$
contains the Weyl group of $\h$. Now we  prove that $N\subset
\Aut(\Delta)$.

 For $x,y\in\gl_m$ put $(x,y)=\tr(xy)$.
The restriction of $(\cdot,\cdot)$ to $\t$ is an $N$-invariant
non-degenerate symmetric bilinear form. Its restriction to
$\t(\mathbb{R})$ is a scalar product.  Further, we notice that the
lattice $X$ generated by the weights of $\varphi$ is invariant
under the action of $N$.

Obviously, $Q\subset X\subset P$.   It follows from
Lemma~\ref{lem3}, that if $X=Q$ and $N\not\subset \Aut(\Delta)$,
then $\Delta=C_4$. Suppose now that $X=P$. Then the dual root
lattice $Q^\vee$ is invariant under the action of $N$ on $\t$.
Lemma~\ref{lem3} implies that if $N\not\subset \Aut(\Delta)$, then
$\Delta=B_4$.

Suppose now that $X\neq Q,P$. Then the group $P/Q$ is not prime.
Tables in~\cite{VO} imply that $\Delta=A_l$, where $l+1$ is not
prime, or $\Delta=D_l$. If $\Delta\neq A_7,A_8,D_8$, then one can
check directly that every element of $P$, whose length is equal to
the length of a root, is a root itself. Therefore if $\Delta\neq
A_7,A_8,D_8$, then $N\subset\Aut(\Delta)$.

The system of weights of the representation $\varphi$ is invariant
under $N\subset \Aut(\Delta)$. Thus, $N$ coincides with the image
of $N_{\widetilde{H}}(\t)$ in $\GL(\t^*)$. So we are done.

Now we construct embeddings of
$\h=\sp_8,\so_9,\sl_8,\sl_9,\so_{16}$, for which $\psi_1$ is not
birational.

 Suppose $\h=\sp_8$.
Let $\varphi:\h\rightarrow \gl_{14}$ be the irreducible
representation with the highest weight $\pi_2$. Let us choose the
orthonormal basis $\varepsilon_1,\ldots, \varepsilon_4\in\t$, so
that $\alpha_i=\varepsilon_i-\varepsilon_{i+1}, i=1,2,3$,
$\alpha_4=2\varepsilon_4$. The weights of $\varphi$ are
$\varepsilon_{i}+\varepsilon_j, i\neq j,$ with multiplicity 1 and
0 with multiplicity 2. The stabilizer of this weight system in
$\GL(\t^*)$ is just $\Aut(D_4)$. Therefore $N_G(\t)/Z_G(\t)\cong
\Aut(D_4)$, while $N_{\overline{H}(\t)}/Z_{\overline{H}}(\t)$ is
the Weyl group of $\Delta$ and has index 3 in $N_G(\t)/Z_G(\t)$.

The algebras $\h=\so_9,\sl_8,\sl_9,\so_{16}$ can be embedded into
the exceptional algebras $\f=F_4,E_7,E_8,E_8$, respectively, as
regular subalgebras. Suppose that $\varphi:\h\hookrightarrow
\gl_n$ is the composition of this embedding and some embedding
$\rho:\f\hookrightarrow \gl_n$. Analogously to the case $\h=\sp_8$
one can show that $N_{\overline{H}}(\t)/Z_{\overline{H}}(\t)$ is
not equal to $N_G(\t)/Z_G(\t)$ because the last group is the Weyl
group of $\f$. This completes the proof of the theorem.
\end{proof}

\section{The algebra $\C[\h]^H$}\label{SECTION_3}
It is known (see~\cite{Bourbaki}, Chapter 8,$\S 8$) that for a
connected semisimple group $H$ the vector space $\C[\h]^H$ is
generated by polynomials $\tr(\rho(x)^n)$, where $\rho$ runs over
the set of all representations of $\h$. In this section we
generalize this result for groups $H$ such that $H^\circ$ is
algebraically simply connected (a reductive algebraic group is
called algebraically simply connected if it is a direct product of
a torus and a simply connected semisimple group).

Let $H$ be a (possibly disconnected)  reductive algebraic group
and $\h$ be its Lie algebra. Denote by $\mathfrak{R}(\h)$ the set
of all equivalence classes of representations of $\h$. Define an
action of the group $H$  on $\mathfrak{R}(\h)$: for $h\in
H,\rho\in\mathfrak{R}(\h)$ we put $(h.\rho)(x)=\rho(\Ad(h)^{-1}x)$
for all $x\in\h$. It is obvious, that
$\mathfrak{R}^{H^\circ}=\mathfrak{R}$. If $\rho$ is the
differential of a representation of  $H$, then $h.\rho=\rho$ for
any $h\in H$.

Suppose $\rho\in\mathfrak{R}(\h)^H$. Then
 $\tr(\rho(x)^n)\in\C[\h]^H$ for all $n\in \N$. Indeed, for $h\in
H,x\in\h$ we have
$h.\tr(\rho(x)^n)=\tr(\rho(\Ad(h)^{-1}x)^n)=\tr((h.\rho)(x)^n)=\tr(\rho(x)^n)$.
\begin{Proposition}\label{app:4.1}
Suppose that  $H^\circ$ is an algebraically simply connected
group. The vector space $\C[\h]^H$ is generated by
$\tr(\rho(x)^n)$, where $\rho$ runs over the set of all
representations of $H$.
\end{Proposition}
\begin{proof}
Suppose the group $H$ is connected. Then $\C[\h]^H\cong
\C[\z(\h)]\otimes \C[\h']^H$, where $\h'=[\h,\h]$. Let $\rho$ be
an irreducible representation of  $\h'$, and $\chi$ be a character
of the torus $Z(H)^\circ$.  Put $\rho_m=\rho\otimes m\chi$. Then
\begin{equation}\label{eq:1}
\tr(\rho_m(x)^n)=\sum_{i=0}^m\begin{pmatrix}n\\i\end{pmatrix}m^{n-i}\chi(x)^{n-i}\tr(\rho(x)^i)
\end{equation}
It follows from (\ref{eq:1}) that  for all $k,l\in\N$ the
polynomial $\chi(x)^k\tr(\rho(x)^l)\in\C[\h]^H$ is a linear
combination of polynomials $\tr(\rho_m(x)^{k+l})$. But the linear
space generated by polynomials of the form
$\chi(x)^k\tr(\rho(x)^l)$, where $\rho,\chi$ are as above,
coincides with $\C[\h]^H$. This completes the proof in the case of
a connected group $H$.

Now we consider the general case. We have $$ \C[\h]^H=\{\sum_{h\in
H/H^\circ} h.f\mid f\in \C[\h]^{H^\circ}\}.$$
 Therefore, the vector space
$\C[\h]^H$ is generated by elements of the form
\begin{equation}\label{eq:2}\sum_{h\in H/H^\circ} h.\tr(\rho(x)^n),
\end{equation}
where $\rho$ is a representation of $H^\circ$. Denote by
$\widetilde{\rho}$ a representation of $H$ induced from $\rho$.
The corresponding representation $\widetilde{\rho}$ of $\h$ is
given by $\widetilde{\rho}=\sum_{h\in H/H^\circ}h.\rho$.   The
polynomial (\ref{eq:2}) is just $\tr(\widetilde{\rho}(x)^n)$.
\end{proof}

Let $I$ be a set of representations of $H$. Denote by $\C[\h^n]^I$
the subalgebra of $\C[\h^n]^H$ generated by polynomials of the
form $\tr(\rho(x)^m)$, where $\rho\in I$. When $I=\{\rho\}$, we
write $\C[\h^n]^\rho$ instead of $\C[\h^n]^{\{\rho\}}$.

It is known from classical invariant theory (see, for example,
\cite{VP}) that if $H=\GL_m,\SL_m,$ $\O_m,\Sp_{2m}$, then
$\C[\h^n]^H=\C[\h^n]^\rho$, where $\rho$ is the tautological
representation of the group $H$. Now let $H$ be one of the
exceptional simple groups $G_2,F_4,E_6,E_7,E_8$ and $\rho$ be the
simplest (=non-trivial irreducible of minimal dimension)
representation of the group $H$. It was shown by several authors
(see~\cite{DQ} for references) that $\C[\h]^H=\C[\h]^\rho$.

\section{Linear equivalence of embeddings}\label{SECTION_4}
Let $G,H$ be  reductive algebraic groups and $\g,\h$ be the Lie
algebras of $G$ and $H$, respectively.  We say that two
homomorphisms $P_1,P_2:H\rightarrow G$ are {\it equivalent} (or,
more precisely, {\it $G$-equivalent}) if there exists $g\in G$
such that $gP_1(x)g^{-1}=P_2(x)$ for all $x\in H^\circ$. Further,
we say that $P_1,P_2$ are {\it linearly equivalent} (or {\it
linearly $G$-equivalent}) if for every representation
$P:G\rightarrow \GL(V)$ representations $P\circ P_1,P\circ P_2$
are $\GL(V)$-equivalent. It is obvious, that equivalent
homomorphisms are linearly equivalent. Homomorphisms
$\rho_1,\rho_2:\h\rightarrow\g$ are said to be $G$-equivalent
(resp., linearly $G$-equivalent) if there exist equivalent (resp.,
linearly equivalent) homomorphisms $P_1,P_2:H\rightarrow G$ with
$dP_1=\rho_1, dP_2=\rho_2$.

Let $P:H\rightarrow G$ be a homomorphism, $\rho:\h\rightarrow \g$
its tangent homomorphism. Denote by $\psi_n^\rho$ the natural
morphism  $\h^n\quo H\rightarrow\rho(\h)^n\quo G$.

The set $\rho(\h)^n\quo G$ can be identified with the set of
equivalence classes $(\rho(x_1),\ldots, \rho(x_n))$, where
$x_1,\ldots, x_n$ generate a reductive subalgebra of $\h$, modulo
$G$-conjugacy. Therefore Lemma~\ref{lem2} implies that
$\psi_n^\rho, n>1,$ is bijective iff for every reductive Lie
algebra $\f$ and embeddings $\rho_1,\rho_2:\f\rightarrow \h$ the
following condition is fulfilled:

if $\rho\circ\rho_1,\rho\circ\rho_2$ are $G$-equivalent, then
$\rho_1,\rho_2$ are $H$-equivalent.

Similarly,  we get the following
\begin{Proposition}\label{app:6.1}
Let $H,G,\h,\g,\rho$ be as above.  The following conditions are
equivalent:
\begin{itemize}
\item[(i)] The map $\psi_1^\rho$ is injective.
\item[(ii)] For every pair $(x_1,x_2)$ of semisimple elements of
 $\h$ if $\rho(x_1)\sim_G\rho(x_2)$, then $x_1\sim_Hx_2$.
\item[(iii)] For any diagonalizable Lie algebra $\t$ and embeddings $\rho_1,\rho_2:\t\hookrightarrow \h$
if $\rho_1,\rho_2$ are $H$-equivalent, then
$\rho\circ\rho_1,\rho\circ\rho_2$ are $G$-equivalent.
\end{itemize}
\end{Proposition}

The following proposition is a generalization of Theorem 1.1.
from~\cite{Dynkin}
\begin{Proposition}\label{app:6.2}
Let $G,\h,\g$ be  as above, $\t\subset\h$ be a Cartan subalgebra,
$P_1,P_2:H\rightarrow G$ be some homomorphisms, $\rho_1=dP_1,
\rho_2=dP_2$.  The following conditions are equivalent
\begin{itemize}
\item[(i)] $\rho_1$ and $\rho_2$ are linearly $G$-equivalent.
\item[(ii)] $\rho_1|_\t$ and $\rho_2|_\t$ are $G$-equivalent.
\end{itemize}
\end{Proposition}
\begin{proof}
$(ii)\Rightarrow (i)$. Replacing $\rho_1$ by $\Ad(g)\circ\rho_1,
g\in G,$ if necessary, we may assume that $\rho_1|_\t=\rho_2|_\t$.
The required result follows from the fact that a representation of
 $\h$ is uniquely determined by the collection of its weights and
their multiplicities.

$(i)\Rightarrow (ii)$. We may assume that $\h=\t$. By
Proposition~\ref{app:6.1}  it is enough to prove the following
claim:

Let $x_1,x_2$ be semisimple elements of  $\g$ such that for every
representation $P:G\rightarrow \GL(V)$ the matrices
$\rho(x_1),\rho(x_2)$ are conjugate, where $\rho=dP$. Then
$x_1,x_2$ are conjugate (with respect to the adjoint action of
$G$).

First we suppose that $G^{\circ}$ is algebraically simply
connected. We see that $\tr(\rho(x_1)^n)=\tr(\rho(x_2)^n)$ for all
$P:G\rightarrow \GL(V)$. It follows from Proposition~\ref{app:4.1}
that $f(x_1)=f(x_2)$ for any $f\in\C[\g]^G$. Since $x_1,x_2$ are
semisimple, we have $x_1\sim_Gx_2$.

Our claim in the general case is now a consequence of the
following
\begin{Lemma}
For every reductive group $G$ there exist a group $\widetilde{G}$
and a surjective covering $\pi:\widetilde{G}\rightarrow G$ such
that $\Ad(\widetilde{G})=\Ad(G)$ and $\widetilde{G}^\circ$ is
algebraically simply connected.
\end{Lemma}
\begin{proof}[Proof of Lemma]
There exist a finite subgroup $F\subset G$ such that
$G=FG^{\circ}$ (see~\cite{Vinberg}, Proposition 7). Let
$\widetilde{G'}$ be the simply connected covering of $(G,G)$. The
group $\widetilde{G}=F\rightthreetimes (Z(G^{\circ})\times
\widetilde{G'})$ and the natural morphism
$\pi:\widetilde{G}\rightarrow G$ has the required properties.
\end{proof}\end{proof}
\begin{Remark}
Let $I$ be a set of representations of $G$ such that the
polynomials $\tr(dP(x)^n), P\in I,$ generate the algebra
$\C[\g]^G$. It follows from the previous proof that
$\rho_1|_\t,\rho_2|_\t$ are $H$-equivalent iff the representations
$dP\circ\rho_1,dP\circ\rho_2$ are $\GL(V)$-equivalent for any
representation $P\in I$.
\end{Remark}

\section{Bijectivity of $\psi_2:\h^2\quo \overline{H}\rightarrow \h^2\quo G$}
\label{SECTION_5}
\begin{Theorem}\label{thm:1}
Let $\h$ be a simple Lie algebra, $H$ an algebraic group with Lie
algebra $\h$. The following conditions are equivalent:
\begin{itemize}
\item[(i)] $\h=\sl_n,\sp_{2n},\so_{2n+1},G_2,F_4$ or
$\h=\so_{2n},n>3,$ and the group $\Ad(H)$ contains an involutory
outer automorphism of the algebra $\h$.
\item[(ii)] For any reductive algebraic Lie algebra $\f$ and any pair of linearly $H$-equivalnent homomorhisms $\rho_1,\rho_2:\f\rightarrow\h$
$\rho_1$ and $\rho_2$ are $H$-equivalent.
\end{itemize}
\end{Theorem}
Theorem~\ref{thm:1} is proved in Sections
\ref{SECTION_7}-\ref{SECTION_10}.
\begin{proof}[Proof of
Theorem~\ref{Thm1:2}]
 It is enough to prove that for every reductive Lie algebra $\f$
 and every embeddings $\rho_1,\rho_2:\f\rightarrow\h$ if
 $\rho_1,\rho_2$ are $G$-equivalent then they are $H$-equivalent.

It follows from the bijectivity of $\psi_1$  that
$\rho_1|_\t,\rho_2|_\t$ are $H$-equivalent. It is enough to prove
that the latter implies $\rho_1,\rho_2$ are $H$-equivalent. One
may assume that $\rho_1|_\t=\rho_2|_\t$. Let
$\h=\z(\h)\oplus\h_1\oplus\ldots\ldots\oplus\h_k$, where $\h_i$ is
a simple noncommutative ideal. Now it's enough to consider the
case $\overline{H}=H_1\times\ldots\times H_k$, where $H_i=\O_{2m}$
if $\h_i=\so_{2m},m>3,$ $H_i=\Int(\h_i)$ otherwise. Denote by
$\pi_i$ the projection from $\h$ to $\h_i$. It's enough to prove
that there is  $h_k\in H_k$ such that
$\Ad(h_k)\circ\pi_k\circ\rho_1=\pi_k\circ\rho_2$. By
Proposition~\ref{app:6.2}, homomorphisms
$\pi_k\circ\rho_1,\pi_k\circ\rho_2$ are linearly $H$-equivalent.
It remains to use Theorem~\ref{thm:1}.
\end{proof}
\begin{Proposition}\label{app:7.1}Let $G$ be a reductive algebraic group such that $G^\circ$ is
simply connected, $\g$ the Lie algebra of $G$. There exists a
representation $P:G\rightarrow \GL(V)$ such that $\rho=dP$ is an
embedding and $\psi_1^\rho:\g\quo G\rightarrow \rho(\g)\quo
\GL(V)$ is injective.
\end{Proposition}
\begin{proof}
We need the following lemma
\begin{Lemma}\label{Lem:1.40}
Let $\rho_1:\g\rightarrow \gl(U_1), \rho_2:\g\rightarrow \gl(U_2)$
be representations, $n=\dim U_1+1$. Put $U=U_1\oplus nU_2$,
$\rho=\rho_1+n\rho_2$. For any two semisimple elements $x,y\in \g$
if the matrices $\rho(x),\rho(y)$ are similar, then for $i=1,2$
the matrices $\rho_i(x),\rho_i(y)$ are similar.
\end{Lemma}
\begin{proof}[Proof of Lemma~\ref{Lem:1.40}]
Denote by $\lambda_1,\ldots,\lambda_k$ the different eigenvalues
of matrices $\rho(x),\rho(y)$ and by $m_1,\ldots,m_k$ their
multiplicities. Both $\rho_1(x)$ and $\rho_1(y)$ (respectively,
$\rho_2(x),\rho_2(y)$) have eigenvalues
$\lambda_1,\ldots,\lambda_k$ with multiplicities
$n\{\frac{m_1}{n}\},\ldots, n\{\frac{m_k}{n}\}$ (respectively,
$[\frac{m_1}{n}],\ldots, [\frac{m_k}{n}]$). Therefore,
$\rho_i(x),\rho_i(y)$ are similar for $i=1,2$.
\end{proof}
It follows from Proposition~\ref{app:4.1} that there exist
representations $P_i:G\rightarrow \GL(U_i), i=\overline{1,m},$
such that the algebra $\C[\g]^G$ is generated by polynomials of
the form $\tr(dP_i(x)^n)$.

For $i=1,\ldots,m$ we define a positive integer $n_i$, a vector
space $\widetilde{U}_i$ and a representation
$\widetilde{P}_i:G\rightarrow \GL(\widetilde{U}_i)$ by formulas
$$n_1=1, \widetilde{U}_1=U_1,\widetilde{P}_1=P_1.$$
$$n_i=\dim\widetilde{U}_{i-1}+1,
\widetilde{U}_i=\widetilde{U}_{i-1}\oplus n_iU_i,
\widetilde{P}_i=\widetilde{P}_{i-1}+n_iP_i.$$

 By Lemma~\ref{Lem:1.40},  for any semisimple $x,y\in \g$ if the matrices
$d\widetilde{P}_m(x),d\widetilde{P}_m(y)$ are similar, then for
$i=\overline{1,m}$ the matrices $dP_i(x),dP_i(y)$ are similar.
Thus $f(x)=f(y)$ for all $f\in \C[\g]^G$ and so $x\sim_G y$. It
follows from Proposition~\ref{app:6.1} that $\psi_1^{dP_m}$ is
bijective.
\end{proof}

\begin{proof}[Proof of Proposition~\ref{Thm1:3}]
First, suppose that $G=\GL_m$. Let $\f$ be a reductive Lie
algebra, $\rho_1,\rho_2:\f\hookrightarrow \h$  linearly
$H$-equivalent but not equivalent embeddings. Such $\rho_1,\rho_2$
exist by Theorem~\ref{thm:1}. If $\psi_2$ is bijective, then
$\rho_1,\rho_2$ are $H$-equivalent. Contradiction.

Now we consider the general case. Since $\h\subset [\g,\g]$, one
may assume that $G$ is a semisimple group. Let $G_1$ be a simply
connected group with Lie algebra $\g$. Changing $G$ with
$G_1\leftthreetimes (\Ad(G)/\Int(\g))$, we may assume that
$G^\circ$ is simply connected.  By Proposition~\ref{app:7.1},
there exists a homomorphism $P:G\rightarrow \GL(V)$ such that
$\rho=dP$ is an embedding and $\psi_1^{\rho}$ is bijective. Since
$\psi_1^\rho$ is injective, it follows that
$\overline{G}:=N_{\GL(V)}(\rho(\g))/Z_{\GL(V)}(\rho(\g))$ is a
group of type I. Therefore, by Theorem~\ref{Thm1:2},
$\psi_2^\rho:\g^2\quo G\rightarrow \g^2\quo \GL(V)$ is bijective.
Denote by $\widetilde{\rho}$ the embedding $\h\hookrightarrow \g$.
Since $\psi_2^{\rho}$ is bijective, the group
$N_{\GL(V)}(\rho(\h))/Z_{\GL(V)}(\rho(\h))\subset \Aut(\h)$
coincides with $\overline{H}$. Now it follows from the first part
of the proof that the map $\psi_2^{\rho\circ\widetilde{\rho}}$ is
not bijective. But
$\psi_2^{\rho\circ\widetilde{\rho}}=\psi_2^\rho\circ\psi_2$.
Therefore $\psi_2$ is not bijective.
\end{proof}

\begin{Proposition}\label{app:6.3}
Let $G,\g,H,\h$ be such as in Proposition~\ref{app:6.1},
$P:H\rightarrow G$ a homomorphism, $\rho=d P$, $\f$ a reductive
Lie algebra, and $\rho_1,\rho_2:\f\rightarrow\h$ embeddings.
Suppose that the equivalent conditions of
Proposition~\ref{app:6.1} are fulfilled  and $\rho_1,\rho_2$ are
$G$-equivalent. If $\f=\sl_2$ or $\rho_1(\f),\rho_2(\f)$ are
regular subalgebras of $\h$, then $\rho_1$ and $\rho_2$ are
$H$-equivalent.
\end{Proposition}
\begin{proof}
Suppose $\f=\sl_2$. Let $(e,h,f)$ be a standard basis of $\f$,
i.e. $[h,e]=2e,[h,f]=-2f,[e,f]=h$. Since
$\rho_1(h)\sim_G\rho_2(h)$ we see that $\rho_1(h)\sim_H\rho_2(h)$.
We may assume that $\rho_1(h)=\rho_2(h)$. There exists $g\in
Z_H(\rho_1(h))$ such that $\Ad(g)\circ\rho_1=\rho_2$ (see, for
example, \cite{Bourbaki}, ch.8, $\S 11$).

Now suppose that $\rho_1(\f)$ and $\rho_2(\f)$ are regular
subalgebras of $\h$. Denote by $\s,\t$  Cartan subalgebras of
$\f,\h$, respectively. There exists $h\in H$ such that
$\Ad(h)\rho_1|_{\s}=\rho_2|_{\s}$. Therefore the proof reduces to
the case when $\rho_1|_\s=\rho_2|_\s$. Since all Cartan
subalgebras of $\z_\g(\rho_1(\s))$ are
$Z_G(\rho_1(\s))$-conjugate, one also may assume that $\t$
normalizes $\rho_1(\f),\rho_2(\f)$. Let $\alpha\in \Delta(\f)$ and
$e_\alpha\in\f^{\alpha}$ be nonzero. Since $\rho_1(\f),\rho_2(\f)$
are regular subalgebras, $\rho^i(e_\alpha)\subset \h^{\beta_i},
i=1,2,$ for some roots $\beta_1,\beta_2$ of $\h$. This implies
$\rho_i(\alpha^\vee)\subset \C\beta_i^\vee$. Since
$\rho^1|_{\s}=\rho^2|_{\s}$, it follows that
$\rho_1(e_{\alpha})=c_\alpha\rho_2(e_{\alpha})$ for some
$c_\alpha\in\C$. There exists $t\in \rho_1(\s)$ such that
$\exp(\ad t)\circ \rho_1=\rho_2$, see~\cite{Bourbaki}, ch.8,$\S
5$.
\end{proof}

\section{Cases $\h=\sl_n,\sp_{2n},\so_{2n+1}$}\label{SECTION_7}
We show that  linearly $H$-equivalent reductive embeddings are
$H$-equivalent for every group $H$ with Lie algebra $\h$.

Let $\f$ be a reductive  Lie algebra and
$\rho_1,\rho_2:\f\hookrightarrow \h$ be linearly $H$-equivalent
embeddings. We have to prove that $\rho_1,\rho_2$ are
$H$-equivalent. By Proposition~\ref{app:6.2}, we may assume that
$\rho_1|_\t=\rho_2|_\t$, where $\t$ is a Cartan subalgebra of
$\f$. Now it is enough to show that $\rho_1,\rho_2$ are
$H^\circ$-equivalent. Denote by $\rho$ the tautological
representation of  $\h$. Recall that
$\C[\h^2]^\rho=\C[\h^2]^{H^\circ}$. It follows  that the
embeddings $\rho_1$ and $\rho_2$ are $H^\circ$-equivalent.

\section{Case $\h=\so_{2n},n>3$}\label{SECTION_8}
First suppose that  $\Ad(H)$ contains an involutory outer
automorphism of  $\h$. Then linearly $H$-equivalent reductive
embeddings are $H$-equivalent. One can prove this analogously to
the previous section using the fact that
$\C[\h^n]^\rho=\C[\h^n]^H$, where $\rho$ is the tautological
representation of $\h$,  $H=\O_{2n}$.

Now suppose that $\Ad(H)=\Int(\h)$. Denote by $\tau$ the
tautological representation of $\h$ and by $\theta$ any outer
involutory automorphism of $\h$. For the proof of the next
proposition see~\cite{Dynkin}.

\begin{Proposition}\label{Prop:9.1}
Let $\f$ be a reductive Lie algebra and $\rho:\f\rightarrow \h$ be
an embedding. Suppose that
\begin{enumerate}
\item The representation $\tau\circ\rho:\f\rightarrow \gl_{2n}$ has zero weight.
\item All irreducible components of $\tau\circ\rho$ have even
dimension.
\end{enumerate}
Then the embeddings $\rho, \theta\circ\rho$ are linearly
$H$-equivalent but not equivalent.
\end{Proposition}

Now it is enough to construct an embedding $\rho:\f\hookrightarrow
\h$ (or representation $\tau\circ\rho$) satisfying the both
conditions of Proposition~\ref{Prop:9.1}. Denote by $\rho_1$ the
adjoint representation of $\sl_3$ (of dimension 8), by $\rho_2$
the exterior square of the tautological representation of $\so_5$
(of dimension 10) and by $\rho_0$ the tautological representation
of $\so_4$. For $n=2k$ we put $\f=\sl_3\oplus\so_4^{k-4}$,
$\tau\circ\rho=\rho_1\oplus (k-4)\rho_0$. For $n=2k+1$ put
$\f=\so_5\oplus\so_4^{k-4},\tau\circ\rho=\rho_2\oplus
(k-4)\rho_0$.

It remains to consider the case  $\h=\so_8, |\Ad(H)/\Int(\h)|=3$.
It is enough to prove that
$\rho_1,\theta\circ\rho_1:\sl_3\hookrightarrow \so_8$ are not
equivalent.

Assume the converse. There exists $h\in H$ such that
$(\Ad(h)\theta)\circ \rho_1=\rho_1$. The order of the image of
$\Ad(h)\theta$ in the group $\Aut(\h)/\Int(\h)$ is 2. This
contradicts Proposition~\ref{Prop:9.1}.

\section{Case $\h=E_l,l=6,7,8$}\label{SECTION_9}
There exists a Levi subalgebra  $\l\subset\h$ isomorphic to
$\so_{10}\times \C^{l-5}$. Put $\f=\so_5\times \C^{l-5}$. Denote
by $\rho^1,\rho^2$ embeddings of $\f$ into $\l$ satisfying the
following conditions
\begin{enumerate}
\item $\rho^1|_{\z(\f)}=\rho^2|_{\z(\f)}$ is an isomorphism of
$\z(\f)$ and $\z(\l)$.
\item $\rho^1|_{\so_5}=\rho_2, \rho^2|_{\so_5}=\theta\circ\rho_2$, where
$\rho_2$ is the exterior square of the tautological representation
of $\so_5$, $\theta$ is an involutory outer automorphism of
$\so_{10}$.
\end{enumerate}

Since $\rho_2,\theta\circ\rho_2:\so_5\hookrightarrow \so_{10}$ are
linearly $\SO_{10}$-equivalent, we see that $\rho^1,\rho^2$ are
linearly $\Int(\h)$-equivalent.

Assume that there exists $h\in \Ad(H)$ such that
$h\circ\rho^1=\rho^2$. Denote by $L'$ a connected subgroup of $H$
with Lie algebra $[\l,\l]$. It is well known that the centralizer
of an algebraic subtorus in a connected reductive algebraic group
is connected. By Proposition~\ref{Prop:9.1}, $h|_{[\l,\l]}\not\in
\Ad(L')$. Since $h\in Z_{H}(\z(\l))$, we have $h\not\in \Int(\h)$.
Thus, $\h=E_6, \Ad(H)=\Aut(\h)$.

Denote by $\t$ a Cartan subalgebra of $\l$. One may assume that
$\t\cap[\l,\l]$ is $\theta$-invariant and that $\theta$ acts on
$\t\cap[\l,\l]$ by a reflection.  The centralizer of
$\rho_2(\so_5)$ in $\SO_{10}$ coincides with the center of
$\SO_{10}$. Hence $\Ad(h)|_{[\l,\l]}=\theta$ and $\Ad(h)|_\t$ is a
reflection because $\Ad(h)$ acts trivially on $\z(\l)$. The
subgroup of $N_H(\t)$ generated by all reflections is the Weyl
group of $\h$. Therefore $h\in\Int(\h)$. Contradiction.

\section{Cases $\h=G_2,F_4$}\label{SECTION_10}
Since the algebra $\h$ has no outer automorphisms, one may assume
that $H$ is connected.

If $\h=G_2$ the statement of Theorem~\ref{thm:1} follows from
Proposition~\ref{app:6.3}. In the sequel we consider the case
$\h=F_4$.

Let $\f$ be a reductive Lie algebra,
$\rho_1,\rho_2:\f\hookrightarrow \h$ linearly $H$-equivalent
embeddings. We have to prove that $\rho_1,\rho_2$ are equivalent.

Assume the converse. It follows from Proposition~\ref{app:6.3}
that $\rank\f<4$.

\begin{Lemma}\label{Lem6}
Let $H$ be a reductive algebraic group,   $\f$  a reductive Lie
algebra such that $\f\cong \s\oplus\f_1$, where $\s,\f_1$ are
ideals of $\f$ and $\rank\s=1$. Suppose
$\rho_1,\rho_2:\f\rightarrow\h$ are linearly $H$-equivalent
embeddings. Then
\begin{itemize}
\item[(1)] There exists $h\in H$ such that $\Ad(h)\circ\rho_1$ and
$\rho_2$ coincide on $\s$.
\item[(2)] Suppose that $\rho_1,\rho_2$ coincide on $\s$. Then
$\rho_1|_{\f_1},\rho_2|_{\f_1}$ are linearly
$Z_H(\rho_1(\s))$-equivalent. \item[(3)] Under assumptions of (2)
$\rho_1,\rho_2$ are $H$-equivalent iff
$\rho_1|_{\f_1},\rho_2|_{\f_1}$ are $Z_H(\rho_1(\s))$-equivalent.
\end{itemize}
\end{Lemma}
\begin{proof}
The third assertion   is obvious. Propositions~\ref{app:6.2}
and~\ref{app:6.3} imply assertion (1),(2) for $\s\cong \C$ and
assertion (1) for $\s=\sl_2$, respectively.

Now one may assume that $\s\cong\sl_2$ and
$\rho_1|_\s=\rho_2|_\s$. Denote by $x$ a general semisimple
element of $\f_1$. By Proposition~\ref{app:6.2},  it is enough to
show that there exists $g\in Z_H(\rho_1(\s))$ such that
$\Ad(g)\rho_1(x)=\rho_2(x)$. Denote by $e,h,f$ a standard basis of
$\rho_1(\s)$. Since $\rho_1,\rho_2$ are linearly $H$-equivalent,
there exist $g_1\in Z_H(h)$ such that
$\Ad(g_1)\rho_1(x)=\rho_2(x)$. Analogously to the proof of
Proposition~\ref{app:6.3}, there exists $g_2\in Z_H(\rho_2(x))$
such that $\Ad(g_2)\Ad(g_1)\rho_1,\rho_2$ coincide on $\s$. The
element $g=g_2g_1$ has the required properties.
\end{proof}

Suppose  that $\f$ is not simple. Then $\f$ satisfies conditions
of  Lemma~\ref{Lem6}. It is easy to see that $Z_H(\rho_1(\s))$ is
a group of type I. Rank of semisimple part of $Z_H(\rho_1(\s))$ is
less than 4. It follows from results of Sections 8,9, that if
$\rho_1,\rho_2$ are linearly equivalent, then they are equivalent
(see the proof of Theorem~\ref{Thm1:2}) . So, we may assume that
$\f$ is simple. Proposition~\ref{app:6.3} implies that
$\rank\f>1$.

Now we prove that the subalgebras $\rho_1(\f),\rho_2(\f)\subset\h$
are not regular. Assume the  converse. To be definite, let
$\rho_1(\f)$ be a regular subalgebra of $\h$. By
Proposition~\ref{app:6.3}, $\rho_2(\f)$ is not regular. Since the
representations $\ad\circ\rho_1,\ad\circ\rho_2$ are equivalent, we
obtain that $\dim\z_\h(\rho_2(\f))=\dim\z_\h(\rho_1(\f))\geqslant
4-\rank\f$. Since for any reductive Lie algebra $\s$ the number
$\dim\s-\rank\s$ is even, $\rank\z_\h(\rho_2(\f))\leqslant
\rank\z_\h(\rho_1(\f))-2$. Therefore, $\rank\z_\h(\rho_2(\f))=0$.
Contradiction.

It follows from~\cite{Dynkin} that every simple subalgebra in $\h$
of rank greater then 1 is contained  in a maximal regular
subalgebra. There are  three maximal regular subalgebras
$\h_1,\h_2,\h_3\subset \h$, $\h_1\cong \so_9, \h_2\cong
\sp_6\times \sl_2, \h_3\cong\sl_3\times \sl_3$. See, for example,
\cite{Dynkin} (there is a mistake in this paper: the subalgebra of
$F_4$ isomorphic to $\sl_4\times \sl_2$ is not maximal, it is
contained in $\h_1$). One may assume that $\h_1,\h_2,\h_3$ contain
a Cartan subalgebra $\t\subset\h$.

All simple subalgebras of rank greater then 1 in $\sp_6$ are
regular. All embeddings $\rho:\f\hookrightarrow \h_1$ such that
the subalgebra $\rho(\f)$ is not regular are listed up to
$\Int(\h_1)$-conjugacy in the following table: \setlongtables
\begin{longtable}{|c|c|}
\hline $\f$&$\rho$\\\hline $\sl_4$&$R(\pi_2)\oplus 3R(0)$\\\hline
$\so_7$&$R(\pi_3)\oplus R(0)$\\\hline $G_2$&$R(\pi_1)\oplus
2R(0)$\\\hline $\so_5$&$R(\pi_2)\oplus R(\pi_2)\oplus
R(0)$\\\hline $\sl_3$&$\ad\oplus R(0)$\\\hline
\end{longtable}
In the second column linear representations in $\C^9$
corresponding to $\rho$ are listed. They determine $\rho$ up to
$\Int(\h_1)$-conjugacy.

Since all simple subalgebras of rank greater then 1 in $\h_3$ are
isomorphic to $\sl_3$, we see that $\f\cong \sl_3$. Denote by
$\rho_1$ the embedding of $\sl_3$ into $\so_9$ listed in the
table. It follows from tables in~\cite{Dynkin} that the
restriction of simplest representation of $\h$ into $\rho_1(\f)$
is isomorphic to $3\ad\oplus 2R(0)$.

It is easy to see that $N_H(\h_3)/N_H(\h_3)^\circ$ is a group of
order 2. The group $N_H(\h_3)$ contains an outer automorphism of
$\h_3$, which acts on $\t$ by multiplication by -1.

Now we introduce some notation. Let $\t$ be a Cartan subalgebra of
$\h$, $\varepsilon_i, i=\overline{1,4}$, its orthonormal basis, so
that $$\Delta(\h)=\{\pm \varepsilon_i\pm\varepsilon_j, i\neq j,
\pm\varepsilon_i, \frac{\pm
\varepsilon_1\pm\varepsilon_2\pm\varepsilon_3\pm\varepsilon_4}{2}
\}.$$ Put $\alpha_1=
(\varepsilon_1-\varepsilon_2-\varepsilon_3-\varepsilon_4)/2,
\alpha_2=\varepsilon_4, \alpha_3=\varepsilon_3-\varepsilon_4,
\alpha_4=\varepsilon_2-\varepsilon_3$. This is a set of simple
roots of $\h$. Denote by $h_1,h_2$ simple coroots of $\f$.

A set of simple roots for $\h_3$ is
$(\varepsilon_1+\varepsilon_2+\varepsilon_3+\varepsilon_4)/2,\alpha_1,\alpha_3,\alpha_4$.
Highest weights for the restriction of the simplest representation
of $\h$ into $\h_3$ are $\varepsilon_1,
(\varepsilon_1+\varepsilon_2-\varepsilon_3-\varepsilon_4)/2,(\varepsilon_1+\varepsilon_2+\varepsilon_3-\varepsilon_4)/2$.

There are two equivalence classes of the embeddings of $\sl_3$
into $\h_3$ up to $N_H(\h_3)$-conjugacy. Only one of these
embeddings is $\GL(26)$-equivalent to $\rho_1$, namely the one
with $h_1\mapsto \varepsilon_1+2\varepsilon_2+\varepsilon_4,\
h_2\mapsto \varepsilon_1-\varepsilon_2-2\varepsilon_4$. We denote
this  embedding by $\rho_2$.

It remains to prove that $\rho_1,\rho_2:\f\hookrightarrow \h$ are
$H$-equivalent. Since an outer automorphism of $\f$ is contained
in both $N_H(\rho_1(\f)),N_H(\rho_2(\f))$, it is enough to show
that the subalgebras $\rho_1(\f),\rho_2(\f)$ are $H$-conjugate.

Denote by $\h_4$ the regular subalgebra of $\h$ corresponding to
the set of roots of maximal length. $\h_4$ is isomorphic to
$\so_8$ and is contained in $\h_1$. One may assume that
$\rho_1(\f)\subset\h_4$. It is clear that $N_H(\t)\subset
N_H(\h_4)$. Therefore $\Ad(N_H(\h_4))=\Aut(\h_4)$. It follows from
description of automorphisms of finite order of simple Lie
algebras (see, for example, \cite{VO}, ch. 4, $\S$4) that there
exists  $t\in N_H(\h_4)$  such that $\Ad_{\h_4}t$ is an element of
order 3 and $\h_4^t=\rho_1(\f)$. $Z_H(\h_4)$ is a finite group of
order 2. Since $t^3\in Z_H(\h_4)$, it follows that $\Ad(t)$ has
finite order (3 or 6). $\h^t$ is a regular subalgebra of rank 4 in
$\h$ and contains $\rho_1(\h)$. Thus $\h^t$ is not contained in a
subalgebra conjugate to $\h_2$.

Let us note that  every proper subalgebra of $\h_1$ containing
$\rho_1(\f)$ coincides with $\h_4$ or $\rho_1(\f)$. Indeed, let
$\widetilde{\f}$ be such a subalgebra. Since centralizer of $\f$
in $\so_9$ is trivial, we obtain that $\widetilde{\f}$ is
semisimple. The algebra $\so_9$ does not contain a subalgebra
isomorphic to $\sl_3\times\sl_3$. Hence, $\f$ is simple. If $\f$
is not conjugate to $\so_8$, then $\f$ is not regular. It follows
from the previous table that $\widetilde{\f}\cong\so_7$ or $\f$. A
subalgebra of $\so_7$ isomorphic to $\sl_3$ is unique up to
conjugation and has non-trivial centralizer. It follows that
$\widetilde{\f}=\f$.

Since order of $t$ is divisible by 3, $\h^t$ is not conjugate to
$\h_4$. Thus, $\h^t$ is not contained in a subalgebra conjugate to
$\h_1$.

  It follows  that $\rho_1(\f)$ is contained in
a subalgebra conjugate to $\h_3$. This completes the proof.

\section{The algebra $\C[\h^n]^{\GL_m}$}
\begin{proof}[Proof of Proposition~\ref{prop:0.1}]
It is known from classical invariant theory that the algebra
$\C[\h^n]^G$ is generated by polynomials $\tr(X_{i_1}X_{i_2}\ldots
X_{i_k})$. Denote by $A$ a subalgebra of the tensor algebra
$T\h^*$ generated by elements of the form
\begin{equation} g(L_1,\ldots,L_d),
\end{equation}
where $g\in S^d(\h^*)^G$, $L_i$ are Lie polynomials in
$X_1,\ldots, X_k$.

Put  $f_k=\tr(X_1\ldots X_k)$. Using polarization we reduce the
required statement to the following one:

 $f_k\in A$ for every positive integer $k$.

The proof of the last statement is by induction on $k$. The case
$k=1$ is trivial. Now assume that we are done for $k<l$.

The symmetric group $S_l$ acts on the space $(\h^*)^{\otimes l}$
by permuting factors in tensor product. It is clear that $A\cap
(\h^*)^{\otimes l}$ is invariant under this action. For every
transposition $\sigma_i=(ii+1)$ one has $$(\sigma_i
f_l)(X_1,\ldots,X_l)-f_l(X_1,\ldots,X_l)=\tr(X_1\ldots
[X_{i},X_{i+1}]\ldots X_k).$$

Therefore, $\sigma_if_l-f_l\in A$. It follows that $f_l-\sigma
f_l\in A$ for every $\sigma\in S_l$. Hence, $f_l\in A$ iff
$$\frac{1}{l!}\sum_{\sigma\in S_l}\sigma f_l\in A.$$ But the
latter is an element of $S^l(\h^*)^G$ and lies in $A$ by
definition.
\end{proof}

\begin{Corollary}
Let $G=\GL_m$, $\h$ be a subalgebra of $\g=\gl_m$,
$\overline{H}=N_G(\h)/Z_G(\h)$. Suppose that
$\psi_1:\C[\h]^G\rightarrow \C[\h]^{\overline{H}}$ is an
isomorphism and  $(\h,\overline{H})$ is one of the following
pairs: $(\sl_k,\Ad(\SL_k))$, $(\so_k,\Ad(\O_k))$,
$(\sp_{2k},\Ad(\Sp_{2k}))$. Then $\psi_n$ is an isomorphism.
\end{Corollary}
\begin{proof}
Let $\rho$ be a representation $\h\hookrightarrow \g$ and $\iota$
be the tautological representation of $\h$.
Proposition~\ref{prop:0.1} implies
$\C[\h^n]^{\rho}=\C[\h^n]^{\iota}$. It follows from classical
invariant theory that $\C[\h^n]^{\iota}=\C[\h^n]^{\overline{H}}$
and we are done.
\end{proof}

\end{document}